\documentclass[11pt,a4paper]{amsart}
\usepackage[T1]{fontenc}
\usepackage[latin9]{inputenc}
\usepackage{amssymb}
\usepackage{color}
\usepackage{graphicx,color}
\usepackage{amsmath, amssymb, graphics}

\makeatletter
\numberwithin{equation}{section} 
\numberwithin{figure}{section} 
  \@ifundefined{theoremstyle}{\usepackage{amsthm}}{}
  \theoremstyle{plain}
  
  \theoremstyle{plain}
  
  \theoremstyle{plain}
  
  \theoremstyle{remark}
  
  \theoremstyle{remark}
  
  \theoremstyle{plain}

\def\com#1{ \hbox{#1}}

\smallskip
\def\<{{\langle }}
\def\>{{\rangle }}

\def\com#1{ \quad\hbox{#1}\quad}

\smallskip
\def\<{{\langle }}
\def\>{{\rangle }}

\begin{document}

\title[Periodic solution Elliptic RTBP]{Near Periodic solution of the Elliptic RTBP for the Jupiter Sun system}
\author{Oscar Perdomo}
\date{\today}

\curraddr{Department of Mathematics\\
Central Connecticut State University\\
New Britain, CT 06050\\}
\email{ perdomoosm@ccsu.edu}


\begin{abstract}

Let us consider the elliptic restricted three body problem (Elliptic RTBP) for the Jupiter Sun system with eccentricity $e=0.048$ and $\mu=0.000953339$. Let us denote by $T$ the period of their orbits. In this paper we provide initial conditions for the position and velocity for a spacecraft such that after one period $T$ the spacecraft comes back to the same place, with   the same velocity,  within an error of 4 meters for the position and 0.2 meters per second for the velocity. Taking this solution as periodic, we present numerical evidence showing that this solution is stable. In order to compare this periodic solution with the motion of celestial bodies in our solar system,  we end this paper by providing an ephemeris of the spacecraft motion from February 17, 2017 to December 28, 2028.
\end{abstract}

\maketitle

\begin{figure}[hbtp]
\begin{center}\includegraphics[width=.3\textwidth]{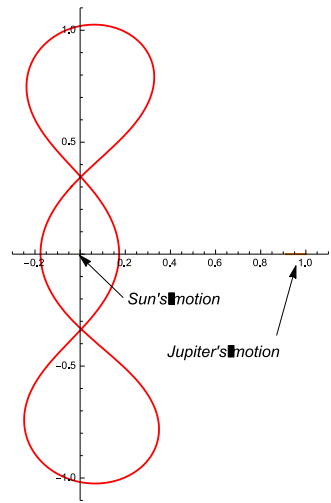}
\end{center}
\caption{Orbit of the spacecraft in a frame that moves along the line that goes from the Sun to Jupiter. In the image, unity is the greatest distance from Jupiter to the Sun, Jupiter's aphelion.
 }\label{fig1}
\end{figure}

\section{Introduction}
As explained in \cite{B1}, finding periodic solutions of the elliptic restricted three body problem is more difficult than finding periodic solutions of the circular restricted three body problem. Among the reasons we have the fact that the elliptic system  is nonconservative and  the fact that the period of a periodic solution in the elliptic case must be a multiple of the period of the motion of the primaries. 

In this paper, rather than extending a solution to a family of solutions or trying to find several families, like it is done, for example in \cite{B}, \cite{B1} and \cite{PYF}, we just take one possible solution and try to obtain the initial conditions that make the trajectory as closed as it can possible be. We notice that even though we made the orbit come back to a position that is within a distance of just 4 meters with a difference in velocity smaller than 0.2 meter per second from its initial velocity, we could not do better after several weeks of computer work. The author is inclined to believe that near his initial conditions there are not initial conditions that make the solution mathematically periodic. This is due to the fact that his approach that started with the regular Newton's method, then with variation of the Newton's method and finally with an exhaustive search, did not seem to converge.

\section{Preliminaries  and main result}

It is well known that for any positive $e<1$, if $\theta(t)$ satisfies the differential equation $\dot{\theta}=\sqrt{\frac{G (m_1+m_2)}{(1-e)^3 r_0^3}}\, (1-e \cos (\theta ))^2$ with $\theta(0)=0$ where $G=6.67408\times 10^{-11}\, {\rm m}^3 {\rm Kg}^{-1} {\rm s}^{-2}$ is the gravitational constant, then, the functions 

$$x(t)=-\frac{(1-e) m_2 r_0}{(m_1+m_2) (1-e \cos (\theta (t)))} \left(\cos (\theta (t)),\sin (\theta (t))\right)$$ 

and 

$$y(t)=\frac{(1-e) m_1 r_0}{(m_1+m_2) (1-e \cos (\theta (t)))}  \left(\cos (\theta (t)),\sin (\theta (t))\right)$$

describe the motion of two bodies with masses $m_1$ and $m_2$ kilograms respectively moving under the effect of the gravity of each other. This is, these two  functions $x$ and $y$ satisfy the system of differential equations

\begin{eqnarray*}
\ddot{x}=\frac{m_2 G}{|y-x|^3} (y-x)\com{and} \ddot{y}=\frac{m_1 G}{|x-y|^3} (x-y)
\end{eqnarray*}

Denoting $\mu=\frac{m_2}{m_1+m_2}$, we have that  $x(0)=-\mu (r_0,0)$, $y(0)= (1-\mu) (r_0,0)$ and the distance between the positions $x$ and $y$ is given by  
$\frac{(1-e) r_0}{1-e \cos (\theta (t))}$. It follows that the maximum distance between the two bodies is $r_0$ meters and it is reached when $t$ equals zero. Also it is known that the period $T$ of the motion of the system equals to $2\pi\sqrt{\frac{r_0^3}{G(1+e)^3\, (m_1+m_2)}}$. 

A direct verification shows that  if we work in the units {\bf ud}, ${\bf ut}$ and ${\bf um}$ such that

$$1{\bf um}=(m1+m_2) {\bf Kg},\quad 1 {\bf ud}=r_0 {\bf m},\quad 1 {\bf ut}=\sqrt{\frac{ r_0^3}{G(m_1+m_2)}} \, {\bf s}$$

then, the period of the motion is $\frac{2\pi}{ (1+e)^\frac{3}{2}} \approx 5.856497259\, $ {\bf ut}, the gravitational constant is $1 {\bf ud}^3 {\bf um}^{-1} {\bf ut}^{-2}$ and, moreover, the expression for $x$ and $y$ reduces to 

$$x(t)=-\frac{(1-e) \mu }{(1-e \cos (\theta (t)))} \left(\cos (\theta (t)),\sin (\theta (t))\right)$$
and

$$y(t)=\frac{(1-e)(1- \mu) }{(1-e \cos (\theta (t)))} \left(\cos (\theta (t)),\sin (\theta (t))\right)$$ 

with $\dot{\theta}=\frac{1}{\sqrt{(1-e)^3}} (1-e \cos (\theta ))^2$. We will represent by $z(t)$ the position of a body with a mass so small that it does not affect the motion of the bodies with masses $m_1$ and $m_2$. Typical examples for the motion described by $z$ are asteroids and spacecraft. If $z(t)$ moves only under the influence of the gravity of the bodies $m_1$ and $m_2$ then, using the new units, it satisfies the following equation,

$$ \ddot{z}=\frac{m_1}{|z-x|^3} (z-x)+ \frac{m_2}{|z-y|^3} (z-y)$$

The main and only result of this paper shows that if 

\begin{eqnarray*}
z(0)&=&(-0.038063861100,0.30182501850)\\
\dot{z}(0)&=&(-1.6227600677,-1.5096541883)
\end{eqnarray*}

then, integrating the differential equation that defines $z$ using the Taylor method of order $9$ with working precision $H=10^{-30}$ and  step $h=0.0005856497259353531319661467$ ($h$ is the 30 digit approximation $\frac{T}{10000}$), we obtain that

\begin{eqnarray*}
z(T)&=&(-0.038063861095882194319990779532, 0.30182501850211801521827874421)\\
\dot{z}(T)&=&(-1.622756783428950379105822092089,-1.5096430093236243034947917510450)
\end{eqnarray*}

Assuming that $r_0=815.757\times 10^{9}$ {\bf m}, $m_1=1.989 \times 10^{30}$ {\bf Kg}, $m_2=1898 \times 10^{24}$ {\bf m}, we obtain that  $1\, {\bf ud}=815.757\times 10^{9}$ {\bf m} and
$1 {\bf ut}=\sqrt{\frac{ r_0^3}{G(m_1+m_2)}} \, {\bf s}=6.39214246027536259802333 \times 10^7$ {\bf s}, and,

\begin{eqnarray*}
z(T)-z(0)&=&(-4.117805680009220468\times 10^{-12},-2.11801521827874421\times 10^{-12}) {\bf ud}\\
&\approx& (-3.35913, -1.72779)\quad \text{\bf meters}
\end{eqnarray*}
{\small
\begin{eqnarray*}
\dot{z}(T)-\dot{z}(0)&=&(-3.284271049620894177907911\times 10^{-6},-0.000011178976375696505208248955) \frac{{\bf ud}}{{\bf ut}}\\
&\approx&(-0.0419135, -0.142665) \quad \text{\bf meters per second}
\end{eqnarray*}
}

We obtain equivalent  results regarding the differences between the initial and final conditions using  the routine from the Software {\it Mathematica} 

{\bf NDSolve[}differential equation,{\bf WorkingPrecision}${\bf \rightarrow 40]}$ 

Treating this motion as a periodic solution, we numerically compute the eigenvalues of the monodromy matrix to be $\lambda_{1}$, $\lambda_{2}$, $\lambda_{3}$ and $\lambda_4$, with $\lambda_{1}$ and $\lambda_{2}$ equal to

$$0.999998796815156697307988946\pm0.001551624627312364108649697 i$$

and $\lambda_{3}$ and $\lambda_4$ equal to 
 
 $$0.974139767581681496571440794\pm 0.225946259107111013549883566 i  $$

We have that 

$$|\lambda_1|=|\lambda_2| =1.0000000005853725651366432083199361400$$

and
$$|\lambda_3|=|\lambda_4| =0.9999999993942950953390829860830566478$$

Since we have four different eigenvalues, assuming they lie in the unit circle, we conclude that the solution is stable. 

\section{Immersing the solution in our solar system from February 17, 2017 to December 28, 2028}

In order to move from the planar solution of the elliptic restricted three body problem presented above to an ephemeris of the solution we need to set up coordinates, then we need to choose an interval of time representing a period for Jupiter and also we need to find the right inclusion $\mathbf{F}$ that sends the two coordinates from the planar solution presented here to the chosen coordinate system for our solar system.

{\bf Period of Jupiter:} We have taken the period of the orbit to be $4332.82$ days. For the planar ODE that we took in the previous section to solve the elliptic restricted three body problem, we have that $t=0$ represents the aphelion of Jupiter. We will be assuming that Jupiter's aphelion happens on February 17, 2017 at 00:00:00.0000. 

{\bf Coordinates for the Solar system :} We will be using one of the solar system used by NASA and we will be using the coordinates displayed by WebHorizon using the following settings:

Center body name: Solar System Barycenter (0)  {source: DE-0431LE-0431}\\
Center-site name: BODY CENTER\\
Center geodetic : 0.00000000,0.00000000,0.0000000 {E-lon(deg),Lat(deg),Alt(km)}\\
Center cylindric: 0.00000000,0.00000000,0.0000000 {E-lon(deg),Dxy(km),Dz(km)}\\
Center radii    : (undefined)   \\                                               
Output units    : AU-D  \\                                                       
Output format   : 02\\
Reference frame : ICRF/J2000.0     \\                                            
Output type     : GEOMETRIC cartesian states\\
Coordinate systm: Ecliptic and Mean Equinox of Reference Epoch  \\      

As an example, the position and velocity of Jupiter on February 17, 2017 is given by:

 {\small A.D. 2017-Feb-17 00:00:00.0000, -5.28409343881439, -1.3382288334361, 0.123732366050597, 0.00176530427744073, \
-0.00695775892189494, -0.0000105642067469229}   

Notice that the information format is: date, position vector in AU, velocity vector in AU/d

{\bf Immersion of the planar coordinates in the solar system coordinates:} In order to find this immersion we found the Jupiter orbital parameters that minimized the difference between the real Jupiter's ephemeris  and the one obtained by assuming that it moves on a perfect ellipse. We have obtained that the best fit is obtained when

\begin{itemize}
\item Jupiter's orbital inclination (deg): 1.300333
\item Jupiter's longitude of ascending node (deg): 100.53
\item Jupiter's argument of perihelion (deg): -86.311
\item Aphelion: 5.453 AU
\end{itemize}

Besides adjusting the orbital parameters, we have  made a translation of 0.007067 AU in the aphelion direction to make the fit better. In this way, if we define

$$\nu=(\cos(100.53-90) \sin(1.30333),\sin(100.53-90)\, \sin(1.30333),\cos(1.30333))$$

$$V_1=(\cos(100.53+180) ,\sin(100.53+180),0)\com{and} V_2=\nu\times V_1$$

and 

$$ {\bf e_1}=\cos(-86.311) V_1+\sin(-86.311)V_2\com{and} {\bf e_2}=\nu\times {\bf e_1} $$

Our immersion is given by $F:R^2\to R^3$ with 

\begin{eqnarray}
F(x_1,x_2)=0.007067 {\bf e_1}+5.453 (x_1 {\bf e_1}+x_2 {\bf e_2}) 
\end{eqnarray}



Recall that $T=\frac{2 \pi}{(1+e)^\frac{3}{2}}\approx 5.856497259\, {\bf ut}$ is the period of the motion in ${\bf ut}$ that we are making equal to $4332.82 $ days. We will be using $F$ to transform positions to the solar system coordinates.

{\bf Verifying the selected immersion:} In order to check the difference between the elliptical motion and the real Jupiter's motion we will consider the following two sequences with 4333 vectors

\begin{eqnarray*}
 H_1&=&\{H_1(1),\dots,H_1(4333)\}\\
 &=& \{F(\frac{(1-e)(1- \mu) }{(1-e \cos (\theta (\frac{i-1}{ UT})))} \left(\cos (\theta (\frac{i-1}{UT})),\sin (\theta (\frac{i-1}{UT}))\right)):\, i=1,2,\dots , 4333\}
 \end{eqnarray*} 

Where $UT=\frac{4332.82}{T}\approx 739.831$, and 

$$ H_2=\{H_2(1),\dots,H_2(4333)\} = \{ (Y_1(t_i),Y_2(t_i),Y_3(t_i)) :  i=1,2,\dots , 4333\}$$

where $(Y_1(t_1),Y_2(t_1),Y_3(t_1))$ are the coordinates of Jupiter on February 17, 2017 at 00:00:00.0000, $ (Y_1(t_2),Y_2(t_2),Y_3(t_2))$ are the coordinates of Jupiter on February 18, 2017 at 00:00:00.0000, and so on. 

A direct verification shows that the distance between the $H_1$ and $H_2$, that is, the sequence of real number given by 

$$dp=\{  |H_1(i)-H_2(i)| , i=1,\dots 4333 \}$$

is a collection of numbers that varies from $0.00237237$ to $0.00444908$ which means the the perfect ellipse in the space given by the immersion $F$ provides an approximation that is within a distance of $0.00444908$ AU for all the days from Feb 17, 2017 to Dec 28, 2028. This approximation is not bad having in mind that everyday Jupiter moves a distance between  $0.007177$  AU and $0.079117$ AU. 

In the same way that $H_1$ generates an ephemeris for Jupiter starting from Feb 17, 2017 which lies on a perfect plane, we can generate an ephemeris for the motion discussed in this paper and, in this way, we can compare this suggested motion for a spacecraft with the motion of other celestial bodies. This ephemeris is displayed on the next section. A video that shows this motion along with the motion of Mercury, Venus, Earth, Mars, Jupiter and the asteroid Hektor is shown in the YouTube video {\color{blue} https://youtu.be/YQ1KY8YmsUA}

\vfil
\eject
\section{Asteroid's ephemeris:2017-Feb-17 00:00:00 to  2017-Feb-17 00:00:00, step: 19 days}
{\tiny
$
\begin{array}{ccccccc}
 \text{ A.D. 2017-Feb-17 00:00:00.0000} & 0.598544 & -1.5462 & -0.00695959 & 6.5468 & 10.1514 & -0.188649 \\
 \text{ A.D. 2017-Mar-08 00:00:00.0000} & 0.758187 & -1.26636 & -0.0116941 & 5.8 & 11.6754 & -0.178281 \\
 \text{ A.D. 2017-Mar-27 00:00:00.0000} & 0.892024 & -0.94467 & -0.0160253 & 4.49126 & 13.3994 & -0.156175 \\
 \text{ A.D. 2017-Apr-15 00:00:00.0000} & 0.98165 & -0.577753 & -0.0195557 & 2.29691 & 15.1449 & -0.114348 \\
 \text{ A.D. 2017-May-04 00:00:00.0000} & 1.0001 & -0.17094 & -0.0216598 & -1.07826 & 16.377 & -0.0439739 \\
 \text{ A.D. 2017-May-23 00:00:00.0000} & 0.918322 & 0.251068 & -0.0215853 & -5.38053 & 16.1678 & 0.0531305 \\
 \text{ A.D. 2017-Jun-11 00:00:00.0000} & 0.726918 & 0.642559 & -0.0189316 & -9.35205 & 14.0279 & 0.150864 \\
 \text{ A.D. 2017-Jun-30 00:00:00.0000} & 0.451992 & 0.962275 & -0.0141113 & -11.7793 & 10.8003 & 0.218579 \\
 \text{ A.D. 2017-Jul-19 00:00:00.0000} & 0.135566 & 1.19854 & -0.00801574 & -12.6695 & 7.67865 & 0.25147 \\
 \text{ A.D. 2017-Aug-07 00:00:00.0000} & -0.190597 & 1.36215 & -0.00140028 & -12.6407 & 5.17 & 0.261257 \\
 \text{ A.D. 2017-Aug-26 00:00:00.0000} & -0.509787 & 1.46949 & 0.00529321 & -12.1876 & 3.27752 & 0.258989 \\
 \text{ A.D. 2017-Sep-14 00:00:00.0000} & -0.814899 & 1.53468 & 0.011847 & -11.5736 & 1.86436 & 0.251131 \\
 \text{ A.D. 2017-Oct-03 00:00:00.0000} & -1.10358 & 1.5683 & 0.0181644 & -10.9206 & 0.798552 & 0.240955 \\
 \text{ A.D. 2017-Oct-22 00:00:00.0000} & -1.37558 & 1.57791 & 0.0242088 & -10.2808 & -0.0184678 & 0.230042 \\
 \text{ A.D. 2017-Nov-10 00:00:00.0000} & -1.63157 & 1.56897 & 0.029972 & -9.67489 & -0.655539 & 0.219137 \\
 \text{ A.D. 2017-Nov-29 00:00:00.0000} & -1.87251 & 1.54546 & 0.0354592 & -9.1092 & -1.1602 & 0.208582 \\
 \text{ A.D. 2017-Dec-18 00:00:00.0000} & -2.09944 & 1.51032 & 0.0406814 & -8.58391 & -1.56557 & 0.198517 \\
 \text{ A.D. 2018-Jan-06 00:00:00.0000} & -2.31338 & 1.46578 & 0.045652 & -8.09663 & -1.89506 & 0.188988 \\
 \text{ A.D. 2018-Jan-25 00:00:00.0000} & -2.51526 & 1.41357 & 0.0503847 & -7.64404 & -2.16557 & 0.179989 \\
 \text{ A.D. 2018-Feb-13 00:00:00.0000} & -2.70592 & 1.35504 & 0.0548929 & -7.22268 & -2.38946 & 0.171494 \\
 \text{ A.D. 2018-Mar-04 00:00:00.0000} & -2.88612 & 1.29125 & 0.059189 & -6.82924 & -2.57599 & 0.163469 \\
 \text{ A.D. 2018-Mar-23 00:00:00.0000} & -3.05655 & 1.22308 & 0.0632846 & -6.46074 & -2.73218 & 0.155876 \\
 \text{ A.D. 2018-Apr-11 00:00:00.0000} & -3.2178 & 1.15122 & 0.0671904 & -6.11449 & -2.86343 & 0.148677 \\
 \text{ A.D. 2018-Apr-30 00:00:00.0000} & -3.37043 & 1.07627 & 0.0709159 & -5.78815 & -2.97394 & 0.141837 \\
 \text{ A.D. 2018-May-19 00:00:00.0000} & -3.5149 & 0.998708 & 0.0744701 & -5.47965 & -3.06706 & 0.135323 \\
 \text{ A.D. 2018-Jun-07 00:00:00.0000} & -3.65166 & 0.918951 & 0.0778608 & -5.1872 & -3.14544 & 0.129107 \\
 \text{ A.D. 2018-Jun-26 00:00:00.0000} & -3.7811 & 0.837346 & 0.0810955 & -4.90919 & -3.21123 & 0.123162 \\
 \text{ A.D. 2018-Jul-15 00:00:00.0000} & -3.90358 & 0.754194 & 0.0841807 & -4.64426 & -3.26619 & 0.117465 \\
 \text{ A.D. 2018-Aug-03 00:00:00.0000} & -4.0194 & 0.669754 & 0.0871225 & -4.39117 & -3.31175 & 0.111993 \\
 \text{ A.D. 2018-Aug-22 00:00:00.0000} & -4.12886 & 0.584252 & 0.0899265 & -4.14885 & -3.34911 & 0.106728 \\
 \text{ A.D. 2018-Sep-10 00:00:00.0000} & -4.23223 & 0.497884 & 0.0925978 & -3.91636 & -3.37925 & 0.101653 \\
 \text{ A.D. 2018-Sep-29 00:00:00.0000} & -4.32974 & 0.410827 & 0.095141 & -3.69283 & -3.40302 & 0.0967517 \\
 \text{ A.D. 2018-Oct-18 00:00:00.0000} & -4.42162 & 0.323232 & 0.0975604 & -3.47752 & -3.42112 & 0.0920109 \\
 \text{ A.D. 2018-Nov-06 00:00:00.0000} & -4.50807 & 0.23524 & 0.09986 & -3.26976 & -3.43415 & 0.0874178 \\
 \text{ A.D. 2018-Nov-25 00:00:00.0000} & -4.58927 & 0.146972 & 0.102043 & -3.06894 & -3.44261 & 0.0829609 \\
 \text{ A.D. 2018-Dec-14 00:00:00.0000} & -4.6654 & 0.0585414 & 0.104114 & -2.87451 & -3.44694 & 0.07863 \\
 \text{ A.D. 2019-Jan-02 00:00:00.0000} & -4.73662 & -0.029952 & 0.106075 & -2.686 & -3.44753 & 0.0744156 \\
 \text{ A.D. 2019-Jan-21 00:00:00.0000} & -4.80306 & -0.118416 & 0.107929 & -2.50294 & -3.44469 & 0.0703092 \\
 \text{ A.D. 2019-Feb-09 00:00:00.0000} & -4.86487 & -0.206766 & 0.109679 & -2.32495 & -3.43871 & 0.0663029 \\
 \text{ A.D. 2019-Feb-28 00:00:00.0000} & -4.92216 & -0.294924 & 0.111327 & -2.15165 & -3.42982 & 0.0623895 \\
 \text{ A.D. 2019-Mar-19 00:00:00.0000} & -4.97507 & -0.38282 & 0.112876 & -1.9827 & -3.41825 & 0.0585623 \\
 \text{ A.D. 2019-Apr-07 00:00:00.0000} & -5.02368 & -0.470386 & 0.114327 & -1.81781 & -3.40417 & 0.0548154 \\
 \text{ A.D. 2019-Apr-26 00:00:00.0000} & -5.06812 & -0.557559 & 0.115684 & -1.65669 & -3.38774 & 0.051143 \\
 \text{ A.D. 2019-May-15 00:00:00.0000} & -5.10845 & -0.644282 & 0.116947 & -1.49907 & -3.36911 & 0.0475399 \\
 \text{ A.D. 2019-Jun-03 00:00:00.0000} & -5.14479 & -0.7305 & 0.118118 & -1.34472 & -3.34839 & 0.0440013 \\
 \text{ A.D. 2019-Jun-22 00:00:00.0000} & -5.1772 & -0.81616 & 0.119199 & -1.19342 & -3.32569 & 0.0405225 \\
 \text{ A.D. 2019-Jul-11 00:00:00.0000} & -5.20576 & -0.901212 & 0.120191 & -1.04496 & -3.3011 & 0.0370995 \\
 \text{ A.D. 2019-Jul-30 00:00:00.0000} & -5.23054 & -0.98561 & 0.121097 & -0.899146 & -3.27471 & 0.0337281 \\
 \text{ A.D. 2019-Aug-18 00:00:00.0000} & -5.25161 & -1.06931 & 0.121916 & -0.755804 & -3.24657 & 0.0304048 \\
 \text{ A.D. 2019-Sep-06 00:00:00.0000} & -5.26903 & -1.15226 & 0.12265 & -0.614765 & -3.21674 & 0.027126 \\
 \text{ A.D. 2019-Sep-25 00:00:00.0000} & -5.28285 & -1.23443 & 0.123301 & -0.475871 & -3.18528 & 0.0238884 \\
 \text{ A.D. 2019-Oct-14 00:00:00.0000} & -5.29314 & -1.31576 & 0.123869 & -0.338976 & -3.15222 & 0.0206888 \\
 \text{ A.D. 2019-Nov-02 00:00:00.0000} & -5.29993 & -1.39623 & 0.124356 & -0.20394 & -3.1176 & 0.0175243 \\
 \text{ A.D. 2019-Nov-21 00:00:00.0000} & -5.30328 & -1.47579 & 0.124762 & -0.0706322 & -3.08144 & 0.0143921 \\
 \text{ A.D. 2019-Dec-10 00:00:00.0000} & -5.30322 & -1.5544 & 0.125087 & 0.0610712 & -3.04375 & 0.0112894 \\
 \text{ A.D. 2019-Dec-29 00:00:00.0000} & -5.2998 & -1.63202 & 0.125333 & 0.191288 & -3.00455 & 0.00821366 \\
 \text{ A.D. 2020-Jan-17 00:00:00.0000} & -5.29305 & -1.70862 & 0.125501 & 0.320132 & -2.96385 & 0.00516242 \\
 \text{ A.D. 2020-Feb-05 00:00:00.0000} & -5.28302 & -1.78416 & 0.125591 & 0.447708 & -2.92165 & 0.00213327 \\
 \text{ A.D. 2020-Feb-24 00:00:00.0000} & -5.26972 & -1.85858 & 0.125603 & 0.574121 & -2.87794 & -0.000876113 \\
 \text{ A.D. 2020-Mar-14 00:00:00.0000} & -5.25318 & -1.93187 & 0.125537 & 0.699467 & -2.83271 & -0.00386797 \\
 \text{ A.D. 2020-Apr-02 00:00:00.0000} & -5.23345 & -2.00398 & 0.125396 & 0.823843 & -2.78595 & -0.00684448 \\
 \text{ A.D. 2020-Apr-21 00:00:00.0000} & -5.21052 & -2.07486 & 0.125178 & 0.947339 & -2.73763 & -0.00980777 \\
 \text{ A.D. 2020-May-10 00:00:00.0000} & -5.18444 & -2.14449 & 0.124884 & 1.07004 & -2.68773 & -0.0127599 \\
 \text{ A.D. 2020-May-29 00:00:00.0000} & -5.15522 & -2.21281 & 0.124514 & 1.19204 & -2.63622 & -0.015703 \\
 \text{ A.D. 2020-Jun-17 00:00:00.0000} & -5.12287 & -2.27979 & 0.124069 & 1.31342 & -2.58307 & -0.018639 \\
 \text{ A.D. 2020-Jul-06 00:00:00.0000} & -5.08741 & -2.34538 & 0.123549 & 1.43425 & -2.52823 & -0.0215699 \\
 \text{ A.D. 2020-Jul-25 00:00:00.0000} & -5.04885 & -2.40954 & 0.122953 & 1.55463 & -2.47165 & -0.0244977 \\
 \text{ A.D. 2020-Aug-13 00:00:00.0000} & -5.00721 & -2.47223 & 0.122282 & 1.67462 & -2.41329 & -0.0274244 \\
 \text{ A.D. 2020-Sep-01 00:00:00.0000} & -4.96249 & -2.53339 & 0.121536 & 1.79431 & -2.35307 & -0.0303519 \\
 \text{ A.D. 2020-Sep-20 00:00:00.0000} & -4.9147 & -2.59299 & 0.120715 & 1.91376 & -2.29095 & -0.0332823 \\
 \text{ A.D. 2020-Oct-09 00:00:00.0000} & -4.86384 & -2.65096 & 0.119818 & 2.03306 & -2.22684 & -0.0362174 \\
 \text{ A.D. 2020-Oct-28 00:00:00.0000} & -4.80992 & -2.70726 & 0.118846 & 2.15228 & -2.16067 & -0.0391593 \\
 \text{ A.D. 2020-Nov-16 00:00:00.0000} & -4.75294 & -2.76183 & 0.117799 & 2.2715 & -2.09234 & -0.04211 \\
 \text{ A.D. 2020-Dec-05 00:00:00.0000} & -4.6929 & -2.81462 & 0.116675 & 2.39078 & -2.02177 & -0.0450716 \\
 \text{ A.D. 2020-Dec-24 00:00:00.0000} & -4.62979 & -2.86556 & 0.115475 & 2.51021 & -1.94885 & -0.0480462 \\
 \text{ A.D. 2021-Jan-12 00:00:00.0000} & -4.56362 & -2.91461 & 0.114199 & 2.62985 & -1.87346 & -0.0510359 \\
 \end{array}$

 \vfil
 \eject
 
 $
 \begin{array}{ccccccc}
 \text{ A.D. 2021-Jan-31 00:00:00.0000} & -4.49436 & -2.96168 & 0.112846 & 2.74979 & -1.79547 & -0.0540429 \\
 \text{ A.D. 2021-Feb-19 00:00:00.0000} & -4.42202 & -3.00671 & 0.111415 & 2.87009 & -1.71476 & -0.0570694 \\
 \text{ A.D. 2021-Mar-10 00:00:00.0000} & -4.34659 & -3.04964 & 0.109906 & 2.99083 & -1.63116 & -0.0601179 \\
 \text{ A.D. 2021-Mar-29 00:00:00.0000} & -4.26805 & -3.09038 & 0.108319 & 3.1121 & -1.54452 & -0.0631907 \\
 \text{ A.D. 2021-Apr-17 00:00:00.0000} & -4.18639 & -3.12885 & 0.106652 & 3.23397 & -1.45465 & -0.0662903 \\
 \text{ A.D. 2021-May-06 00:00:00.0000} & -4.10158 & -3.16498 & 0.104905 & 3.35651 & -1.36136 & -0.0694193 \\
 \text{ A.D. 2021-May-25 00:00:00.0000} & -4.01363 & -3.19866 & 0.103078 & 3.47982 & -1.26442 & -0.0725805 \\
 \text{ A.D. 2021-Jun-13 00:00:00.0000} & -3.92249 & -3.2298 & 0.101169 & 3.60396 & -1.1636 & -0.0757765 \\
 \text{ A.D. 2021-Jul-02 00:00:00.0000} & -3.82816 & -3.2583 & 0.0991771 & 3.72902 & -1.05862 & -0.0790105 \\
 \text{ A.D. 2021-Jul-21 00:00:00.0000} & -3.7306 & -3.28405 & 0.0971019 & 3.85509 & -0.9492 & -0.0822854 \\
 \text{ A.D. 2021-Aug-09 00:00:00.0000} & -3.62979 & -3.30692 & 0.094942 & 3.98225 & -0.835006 & -0.0856046 \\
 \text{ A.D. 2021-Aug-28 00:00:00.0000} & -3.5257 & -3.3268 & 0.0926963 & 4.11059 & -0.715677 & -0.0889713 \\
 \text{ A.D. 2021-Sep-16 00:00:00.0000} & -3.41829 & -3.34355 & 0.0903635 & 4.24018 & -0.59081 & -0.0923893 \\
 \text{ A.D. 2021-Oct-05 00:00:00.0000} & -3.30754 & -3.35701 & 0.0879422 & 4.37112 & -0.459951 & -0.0958623 \\
 \text{ A.D. 2021-Oct-24 00:00:00.0000} & -3.19341 & -3.36703 & 0.085431 & 4.50348 & -0.322595 & -0.0993942 \\
 \text{ A.D. 2021-Nov-12 00:00:00.0000} & -3.07587 & -3.37343 & 0.0828282 & 4.63735 & -0.178172 & -0.102989 \\
 \text{ A.D. 2021-Dec-01 00:00:00.0000} & -2.95486 & -3.37602 & 0.0801323 & 4.77281 & -0.0260392 & -0.106651 \\
 \text{ A.D. 2021-Dec-20 00:00:00.0000} & -2.83035 & -3.3746 & 0.0773414 & 4.9099 & 0.134533 & -0.110386 \\
 \text{ A.D. 2022-Jan-08 00:00:00.0000} & -2.70231 & -3.36894 & 0.0744537 & 5.04871 & 0.304375 & -0.114197 \\
 \text{ A.D. 2022-Jan-27 00:00:00.0000} & -2.57067 & -3.35879 & 0.071467 & 5.18925 & 0.484435 & -0.118089 \\
 \text{ A.D. 2022-Feb-15 00:00:00.0000} & -2.4354 & -3.34388 & 0.0683793 & 5.33156 & 0.675803 & -0.122068 \\
 \text{ A.D. 2022-Mar-06 00:00:00.0000} & -2.29646 & -3.32388 & 0.0651882 & 5.47561 & 0.879737 & -0.126138 \\
 \text{ A.D. 2022-Mar-25 00:00:00.0000} & -2.15379 & -3.29848 & 0.0618913 & 5.62136 & 1.0977 & -0.130305 \\
 \text{ A.D. 2022-Apr-13 00:00:00.0000} & -2.00736 & -3.26728 & 0.0584862 & 5.76868 & 1.33139 & -0.134571 \\
 \text{ A.D. 2022-May-02 00:00:00.0000} & -1.85713 & -3.22985 & 0.0549702 & 5.91736 & 1.58282 & -0.138943 \\
 \text{ A.D. 2022-May-21 00:00:00.0000} & -1.70307 & -3.18572 & 0.0513406 & 6.06709 & 1.85436 & -0.143421 \\
 \text{ A.D. 2022-Jun-09 00:00:00.0000} & -1.54515 & -3.13432 & 0.0475945 & 6.21737 & 2.14883 & -0.148007 \\
 \text{ A.D. 2022-Jun-28 00:00:00.0000} & -1.38338 & -3.07504 & 0.0437293 & 6.36748 & 2.46963 & -0.152698 \\
 \text{ A.D. 2022-Jul-17 00:00:00.0000} & -1.21776 & -3.00713 & 0.0397424 & 6.51633 & 2.82087 & -0.157488 \\
 \text{ A.D. 2022-Aug-05 00:00:00.0000} & -1.04835 & -2.92975 & 0.0356313 & 6.66238 & 3.20758 & -0.162363 \\
 \text{ A.D. 2022-Aug-24 00:00:00.0000} & -0.875248 & -2.84193 & 0.0313942 & 6.80334 & 3.63594 & -0.167297 \\
 \text{ A.D. 2022-Sep-12 00:00:00.0000} & -0.698628 & -2.74249 & 0.02703 & 6.93585 & 4.11369 & -0.172248 \\
 \text{ A.D. 2022-Oct-01 00:00:00.0000} & -0.518764 & -2.63005 & 0.0225392 & 7.05493 & 4.65051 & -0.177143 \\
 \text{ A.D. 2022-Oct-20 00:00:00.0000} & -0.336091 & -2.50293 & 0.0179246 & 7.15307 & 5.2587 & -0.181867 \\
 \text{ A.D. 2022-Nov-08 00:00:00.0000} & -0.151284 & -2.35911 & 0.0131928 & 7.21881 & 5.95397 & -0.186229 \\
 \text{ A.D. 2022-Nov-27 00:00:00.0000} & 0.034613 & -2.19611 & 0.00835684 & 7.23424 & 6.7564 & -0.18991 \\
 \text{ A.D. 2022-Dec-16 00:00:00.0000} & 0.21997 & -2.01086 & 0.00344046 & 7.17078 & 7.69147 & -0.192378 \\
 \text{ A.D. 2023-Jan-04 00:00:00.0000} & 0.402211 & -1.79956 & -0.00151451 & 6.98167 & 8.79048 & -0.192718 \\
 \text{ A.D. 2023-Jan-23 00:00:00.0000} & 0.577192 & -1.55756 & -0.00643475 & 6.58874 & 10.0882 & -0.189325 \\
 \text{ A.D. 2023-Feb-11 00:00:00.0000} & 0.738128 & -1.27936 & -0.0111913 & 5.85964 & 11.6117 & -0.17935 \\
 \text{ A.D. 2023-Mar-02 00:00:00.0000} & 0.873788 & -0.959249 & -0.0155568 & 4.5746 & 13.3418 & -0.157799 \\
 \text{ A.D. 2023-Mar-21 00:00:00.0000} & 0.965902 & -0.593588 & -0.0191376 & 2.40719 & 15.1079 & -0.116661 \\
 \text{ A.D. 2023-Apr-09 00:00:00.0000} & 0.987436 & -0.187222 & -0.0213089 & -0.952323 & 16.383 & -0.046816 \\
 \text{ A.D. 2023-Apr-28 00:00:00.0000} & 0.908712 & 0.235637 & -0.0213061 & -5.27666 & 16.2265 & 0.0505631 \\
 \text{ A.D. 2023-May-17 00:00:00.0000} & 0.71927 & 0.62901 & -0.0187042 & -9.30561 & 14.1077 & 0.149493 \\
 \text{ A.D. 2023-Jun-05 00:00:00.0000} & 0.44486 & 0.95058 & -0.0139032 & -11.7808 & 10.8603 & 0.218362 \\
 \text{ A.D. 2023-Jun-24 00:00:00.0000} & 0.128099 & 1.188 & -0.00780486 & -12.69 & 7.70858 & 0.251804 \\
 \text{ A.D. 2023-Jul-13 00:00:00.0000} & -0.198634 & 1.35207 & -0.00117858 & -12.6629 & 5.17779 & 0.261719 \\
 \text{ A.D. 2023-Aug-01 00:00:00.0000} & -0.51834 & 1.45942 & 0.00552638 & -12.2052 & 3.27228 & 0.259404 \\
 \text{ A.D. 2023-Aug-20 00:00:00.0000} & -0.823832 & 1.52438 & 0.0120896 & -11.5856 & 1.85217 & 0.251451 \\
 \text{ A.D. 2023-Sep-08 00:00:00.0000} & -1.11275 & 1.55763 & 0.018414 & -10.9276 & 0.782924 & 0.241177 \\
 \text{ A.D. 2023-Sep-27 00:00:00.0000} & -1.38488 & 1.56682 & 0.024463 & -10.2836 & -0.0355688 & 0.230176 \\
 \text{ A.D. 2023-Oct-16 00:00:00.0000} & -1.6409 & 1.55744 & 0.0302286 & -9.67438 & -0.673026 & 0.219198 \\
 \text{ A.D. 2023-Nov-04 00:00:00.0000} & -1.88179 & 1.53348 & 0.0357166 & -9.10602 & -1.17747 & 0.208582 \\
 \text{ A.D. 2023-Nov-23 00:00:00.0000} & -2.10861 & 1.4979 & 0.0409381 & -8.57861 & -1.58229 & 0.198468 \\
 \text{ A.D. 2023-Dec-12 00:00:00.0000} & -2.32239 & 1.45294 & 0.0459069 & -8.08962 & -1.91107 & 0.188897 \\
 \text{ A.D. 2023-Dec-31 00:00:00.0000} & -2.52407 & 1.40033 & 0.0506369 & -7.63565 & -2.18076 & 0.179864 \\
 \text{ A.D. 2024-Jan-19 00:00:00.0000} & -2.7145 & 1.34142 & 0.0551414 & -7.21315 & -2.40381 & 0.171341 \\
 \text{ A.D. 2024-Feb-07 00:00:00.0000} & -2.89445 & 1.27727 & 0.0594333 & -6.81878 & -2.5895 & 0.163292 \\
 \text{ A.D. 2024-Feb-26 00:00:00.0000} & -3.06459 & 1.20876 & 0.0635241 & -6.4495 & -2.74484 & 0.155677 \\
 \text{ A.D. 2024-Mar-16 00:00:00.0000} & -3.22555 & 1.13659 & 0.0674245 & -6.10261 & -2.87527 & 0.14846 \\
 \text{ A.D. 2024-Apr-04 00:00:00.0000} & -3.37786 & 1.06135 & 0.0711443 & -5.77574 & -2.98499 & 0.141605 \\
 \text{ A.D. 2024-Apr-23 00:00:00.0000} & -3.52201 & 0.983513 & 0.0746923 & -5.46679 & -3.07734 & 0.135078 \\
 \text{ A.D. 2024-May-12 00:00:00.0000} & -3.65844 & 0.9035 & 0.0780766 & -5.17395 & -3.15499 & 0.128851 \\
 \text{ A.D. 2024-May-31 00:00:00.0000} & -3.78754 & 0.821659 & 0.0813046 & -4.89564 & -3.22008 & 0.122896 \\
 \text{ A.D. 2024-Jun-19 00:00:00.0000} & -3.90966 & 0.738289 & 0.0843828 & -4.63044 & -3.27436 & 0.11719 \\
 \text{ A.D. 2024-Jul-08 00:00:00.0000} & -4.02512 & 0.653647 & 0.0873175 & -4.37713 & -3.31927 & 0.11171 \\
 \text{ A.D. 2024-Jul-27 00:00:00.0000} & -4.13422 & 0.56796 & 0.0901142 & -4.13463 & -3.356 & 0.106439 \\
 \text{ A.D. 2024-Aug-15 00:00:00.0000} & -4.23722 & 0.481423 & 0.0927779 & -3.90198 & -3.38555 & 0.101358 \\
 \text{ A.D. 2024-Sep-03 00:00:00.0000} & -4.33436 & 0.394211 & 0.0953134 & -3.67834 & -3.40875 & 0.0964514 \\
 \text{ A.D. 2024-Sep-22 00:00:00.0000} & -4.42587 & 0.306477 & 0.0977251 & -3.46293 & -3.4263 & 0.0917061 \\
 \text{ A.D. 2024-Oct-11 00:00:00.0000} & -4.51194 & 0.218358 & 0.100017 & -3.25509 & -3.43879 & 0.087109 \\
 \text{ A.D. 2024-Oct-30 00:00:00.0000} & -4.59277 & 0.129978 & 0.102192 & -3.05421 & -3.44675 & 0.0826488 \\
 \text{ A.D. 2024-Nov-18 00:00:00.0000} & -4.66852 & 0.0414471 & 0.104255 & -2.85975 & -3.45059 & 0.0783149 \\
 \text{ A.D. 2024-Dec-07 00:00:00.0000} & -4.73935 & -0.047134 & 0.106207 & -2.67121 & -3.45071 & 0.074098 \\
 \text{ A.D. 2024-Dec-26 00:00:00.0000} & -4.80541 & -0.135674 & 0.108053 & -2.48814 & -3.44741 & 0.0699894 \\
 \text{ A.D. 2025-Jan-14 00:00:00.0000} & -4.86684 & -0.224088 & 0.109795 & -2.31014 & -3.44099 & 0.0659812 \\
 \text{ A.D. 2025-Feb-02 00:00:00.0000} & -4.92376 & -0.3123 & 0.111435 & -2.13685 & -3.43168 & 0.0620662 \\
 \text{ A.D. 2025-Feb-21 00:00:00.0000} & -4.97628 & -0.400237 & 0.112975 & -1.96792 & -3.4197 & 0.0582377 \\
 \text{ A.D. 2025-Mar-12 00:00:00.0000} & -5.02452 & -0.487835 & 0.114419 & -1.80305 & -3.40522 & 0.0544897 \\
 \text{ A.D. 2025-Mar-31 00:00:00.0000} & -5.06857 & -0.575031 & 0.115767 & -1.64196 & -3.3884 & 0.0508164 \\
 \text{ A.D. 2025-Apr-19 00:00:00.0000} & -5.10854 & -0.661766 & 0.117021 & -1.48438 & -3.3694 & 0.0472126 \\
 \end{array}$
 
 \vfil
 \eject
 
 $\begin{array}{ccccccc}
 \text{ A.D. 2025-May-08 00:00:00.0000} & -5.14449 & -0.747986 & 0.118184 & -1.33008 & -3.34832 & 0.0436734 \\
 \text{ A.D. 2025-May-27 00:00:00.0000} & -5.17653 & -0.83364 & 0.119256 & -1.17883 & -3.32526 & 0.0401943 \\
 \text{ A.D. 2025-Jun-15 00:00:00.0000} & -5.20471 & -0.918677 & 0.120241 & -1.03042 & -3.30033 & 0.036771 \\
 \text{ A.D. 2025-Jul-04 00:00:00.0000} & -5.22912 & -1.00305 & 0.121137 & -0.884664 & -3.2736 & 0.0333996 \\
 \text{ A.D. 2025-Jul-23 00:00:00.0000} & -5.24982 & -1.08671 & 0.121948 & -0.741384 & -3.24513 & 0.0300763 \\
 \text{ A.D. 2025-Aug-11 00:00:00.0000} & -5.26687 & -1.16963 & 0.122674 & -0.60041 & -3.21499 & 0.0267976 \\
 \text{ A.D. 2025-Aug-30 00:00:00.0000} & -5.28033 & -1.25174 & 0.123317 & -0.461586 & -3.18321 & 0.0235602 \\
 \text{ A.D. 2025-Sep-18 00:00:00.0000} & -5.29024 & -1.33302 & 0.123877 & -0.324763 & -3.14985 & 0.020361 \\
 \text{ A.D. 2025-Oct-07 00:00:00.0000} & -5.29667 & -1.41343 & 0.124355 & -0.189802 & -3.11492 & 0.0171969 \\
 \text{ A.D. 2025-Oct-26 00:00:00.0000} & -5.29966 & -1.49291 & 0.124752 & -0.0565731 & -3.07846 & 0.0140652 \\
 \text{ A.D. 2025-Nov-14 00:00:00.0000} & -5.29924 & -1.57144 & 0.125069 & 0.075049 & -3.04048 & 0.0109631 \\
 \text{ A.D. 2025-Dec-03 00:00:00.0000} & -5.29546 & -1.64898 & 0.125307 & 0.205182 & -3.001 & 0.00788811 \\
 \text{ A.D. 2025-Dec-22 00:00:00.0000} & -5.28836 & -1.72548 & 0.125466 & 0.333939 & -2.96002 & 0.00483763 \\
 \text{ A.D. 2026-Jan-10 00:00:00.0000} & -5.27797 & -1.80091 & 0.125547 & 0.461426 & -2.91754 & 0.00180931 \\
 \text{ A.D. 2026-Jan-29 00:00:00.0000} & -5.26432 & -1.87523 & 0.125551 & 0.587746 & -2.87355 & -0.00119916 \\
 \text{ A.D. 2026-Feb-17 00:00:00.0000} & -5.24744 & -1.94841 & 0.125478 & 0.712998 & -2.82804 & -0.00419003 \\
 \text{ A.D. 2026-Mar-08 00:00:00.0000} & -5.22735 & -2.02039 & 0.125328 & 0.837276 & -2.78101 & -0.00716548 \\
 \text{ A.D. 2026-Mar-27 00:00:00.0000} & -5.20409 & -2.09114 & 0.125102 & 0.960672 & -2.73242 & -0.0101276 \\
 \text{ A.D. 2026-Apr-15 00:00:00.0000} & -5.17767 & -2.16063 & 0.124799 & 1.08327 & -2.68226 & -0.0130786 \\
 \text{ A.D. 2026-May-04 00:00:00.0000} & -5.1481 & -2.22881 & 0.124422 & 1.20516 & -2.63049 & -0.0160204 \\
 \text{ A.D. 2026-May-23 00:00:00.0000} & -5.11542 & -2.29564 & 0.123968 & 1.32643 & -2.57707 & -0.018955 \\
 \text{ A.D. 2026-Jun-11 00:00:00.0000} & -5.07963 & -2.36107 & 0.12344 & 1.44715 & -2.52196 & -0.0218845 \\
 \text{ A.D. 2026-Jun-30 00:00:00.0000} & -5.04074 & -2.42507 & 0.122836 & 1.56741 & -2.46512 & -0.0248109 \\
 \text{ A.D. 2026-Jul-19 00:00:00.0000} & -4.99877 & -2.48758 & 0.122157 & 1.68728 & -2.40649 & -0.0277359 \\
 \text{ A.D. 2026-Aug-07 00:00:00.0000} & -4.95373 & -2.54857 & 0.121403 & 1.80685 & -2.34602 & -0.0306618 \\
 \text{ A.D. 2026-Aug-26 00:00:00.0000} & -4.90562 & -2.60798 & 0.120574 & 1.92617 & -2.28363 & -0.0335904 \\
 \text{ A.D. 2026-Sep-14 00:00:00.0000} & -4.85444 & -2.66576 & 0.11967 & 2.04534 & -2.21926 & -0.0365236 \\
 \text{ A.D. 2026-Oct-03 00:00:00.0000} & -4.80021 & -2.72186 & 0.11869 & 2.16443 & -2.15282 & -0.0394636 \\
 \text{ A.D. 2026-Oct-22 00:00:00.0000} & -4.74292 & -2.77623 & 0.117634 & 2.2835 & -2.08422 & -0.0424122 \\
 \text{ A.D. 2026-Nov-10 00:00:00.0000} & -4.68257 & -2.8288 & 0.116503 & 2.40264 & -2.01338 & -0.0453717 \\
 \text{ A.D. 2026-Nov-29 00:00:00.0000} & -4.61916 & -2.87953 & 0.115296 & 2.52191 & -1.94019 & -0.0483439 \\
 \text{ A.D. 2026-Dec-18 00:00:00.0000} & -4.55268 & -2.92835 & 0.114012 & 2.64139 & -1.86452 & -0.0513312 \\
 \text{ A.D. 2027-Jan-06 00:00:00.0000} & -4.48314 & -2.97519 & 0.112651 & 2.76116 & -1.78626 & -0.0543356 \\
 \text{ A.D. 2027-Jan-25 00:00:00.0000} & -4.41051 & -3.01998 & 0.111212 & 2.88129 & -1.70527 & -0.0573595 \\
 \text{ A.D. 2027-Feb-13 00:00:00.0000} & -4.33479 & -3.06266 & 0.109696 & 3.00186 & -1.6214 & -0.0604051 \\
 \text{ A.D. 2027-Mar-04 00:00:00.0000} & -4.25597 & -3.10315 & 0.108101 & 3.12294 & -1.53447 & -0.0634749 \\
 \text{ A.D. 2027-Mar-23 00:00:00.0000} & -4.17403 & -3.14136 & 0.106427 & 3.24461 & -1.44432 & -0.0665712 \\
 \text{ A.D. 2027-Apr-11 00:00:00.0000} & -4.08896 & -3.17721 & 0.104674 & 3.36694 & -1.35073 & -0.0696968 \\
 \text{ A.D. 2027-Apr-30 00:00:00.0000} & -4.00073 & -3.21062 & 0.102839 & 3.49003 & -1.2535 & -0.0728543 \\
 \text{ A.D. 2027-May-19 00:00:00.0000} & -3.90934 & -3.24147 & 0.100923 & 3.61394 & -1.15237 & -0.0760465 \\
 \text{ A.D. 2027-Jun-07 00:00:00.0000} & -3.81475 & -3.26968 & 0.0989245 & 3.73876 & -1.04709 & -0.0792763 \\
 \text{ A.D. 2027-Jun-26 00:00:00.0000} & -3.71694 & -3.29513 & 0.0968425 & 3.86458 & -0.937364 & -0.0825468 \\
 \text{ A.D. 2027-Jul-15 00:00:00.0000} & -3.61589 & -3.3177 & 0.094676 & 3.99147 & -0.822857 & -0.0858613 \\
 \text{ A.D. 2027-Aug-03 00:00:00.0000} & -3.51157 & -3.33726 & 0.0924238 & 4.11952 & -0.703211 & -0.089223 \\
 \text{ A.D. 2027-Aug-22 00:00:00.0000} & -3.40394 & -3.35368 & 0.0900846 & 4.24881 & -0.578021 & -0.0926356 \\
 \text{ A.D. 2027-Sep-10 00:00:00.0000} & -3.29297 & -3.36681 & 0.087657 & 4.37943 & -0.446836 & -0.0961027 \\
 \text{ A.D. 2027-Sep-29 00:00:00.0000} & -3.17863 & -3.37649 & 0.0851397 & 4.51145 & -0.309149 & -0.0996284 \\
 \text{ A.D. 2027-Oct-18 00:00:00.0000} & -3.06089 & -3.38254 & 0.082531 & 4.64496 & -0.16439 & -0.103217 \\
 \text{ A.D. 2027-Nov-06 00:00:00.0000} & -2.93969 & -3.38477 & 0.0798294 & 4.78002 & -0.0119172 & -0.106872 \\
 \text{ A.D. 2027-Nov-25 00:00:00.0000} & -2.815 & -3.38299 & 0.0770329 & 4.9167 & 0.148999 & -0.110598 \\
 \text{ A.D. 2027-Dec-14 00:00:00.0000} & -2.68679 & -3.37695 & 0.0741398 & 5.05506 & 0.319187 & -0.1144 \\
 \text{ A.D. 2028-Jan-02 00:00:00.0000} & -2.55499 & -3.36642 & 0.071148 & 5.19513 & 0.499595 & -0.118284 \\
 \text{ A.D. 2028-Jan-21 00:00:00.0000} & -2.41958 & -3.35111 & 0.0680554 & 5.33692 & 0.691312 & -0.122252 \\
 \text{ A.D. 2028-Feb-09 00:00:00.0000} & -2.28051 & -3.33071 & 0.0648598 & 5.48041 & 0.895594 & -0.126312 \\
 \text{ A.D. 2028-Feb-28 00:00:00.0000} & -2.13773 & -3.3049 & 0.0615586 & 5.62555 & 1.1139 & -0.130466 \\
 \text{ A.D. 2028-Mar-18 00:00:00.0000} & -1.9912 & -3.27328 & 0.0581496 & 5.77221 & 1.34793 & -0.134719 \\
 \text{ A.D. 2028-Apr-06 00:00:00.0000} & -1.84088 & -3.23542 & 0.0546299 & 5.92018 & 1.59969 & -0.139076 \\
 \text{ A.D. 2028-Apr-25 00:00:00.0000} & -1.68676 & -3.19085 & 0.0509971 & 6.06913 & 1.87154 & -0.143538 \\
 \text{ A.D. 2028-May-14 00:00:00.0000} & -1.5288 & -3.13901 & 0.0472482 & 6.21855 & 2.16631 & -0.148106 \\
 \text{ A.D. 2028-Jun-02 00:00:00.0000} & -1.367 & -3.07927 & 0.0433808 & 6.36771 & 2.48737 & -0.152777 \\
 \text{ A.D. 2028-Jun-21 00:00:00.0000} & -1.20139 & -3.0109 & 0.0393921 & 6.51552 & 2.83883 & -0.157545 \\
 \text{ A.D. 2028-Jul-10 00:00:00.0000} & -1.03202 & -2.93306 & 0.0352798 & 6.66041 & 3.22569 & -0.162394 \\
 \text{ A.D. 2028-Jul-29 00:00:00.0000} & -0.858988 & -2.84477 & 0.0310423 & 6.80007 & 3.65412 & -0.167299 \\
 \text{ A.D. 2028-Aug-17 00:00:00.0000} & -0.68247 & -2.74487 & 0.0266785 & 6.93112 & 4.1318 & -0.172217 \\
 \text{ A.D. 2028-Sep-05 00:00:00.0000} & -0.502749 & -2.63196 & 0.0221889 & 7.04856 & 4.66837 & -0.177075 \\
 \text{ A.D. 2028-Sep-24 00:00:00.0000} & -0.320263 & -2.50439 & 0.0175766 & 7.14487 & 5.27603 & -0.181756 \\
 \text{ A.D. 2028-Oct-13 00:00:00.0000} & -0.135692 & -2.36014 & 0.0128483 & 7.20855 & 5.97036 & -0.186067 \\
 \text{ A.D. 2028-Nov-01 00:00:00.0000} & 0.049913 & -2.19673 & 0.0080172 & 7.22173 & 6.77124 & -0.189692 \\
 \text{ A.D. 2028-Nov-20 00:00:00.0000} & 0.234918 & -2.01113 & 0.00310723 & 7.15589 & 7.70383 & -0.192097 \\
 \text{ A.D. 2028-Dec-09 00:00:00.0000} & 0.416747 & -1.79956 & -0.00183966 & 6.96454 & 8.79894 & -0.19237 \\
 \text{ A.D. 2028-Dec-28 00:00:00.0000} & 0.591267 & -1.55742 & -0.00675017 & 6.57014 & 10.0906 & -0.188918 \\
\end{array}
$
}




\begin{thebibliography}{11}


\bibitem{B} R. Broucke.  \emph{Stability of Periodic Orbits in the Elliptic, Restricted Tree-Body Problem}. AIAA Journal, Vol 7, No. 6, June 1969. pp. 1003-1009.


\bibitem{B1} R. Broucke.  \emph{On the Elliptic Restricted Three-Body Problem}. The Journal of the Astronautical Sciences, Vol XIX, No. 6, pp. 417-432.


\bibitem{PYF} J.F. Palacian, P. Yanguas, S.  Fern\'andez.  \emph{Searching for periodic orbits of the spatial elliptic restricted three body problem by double averaging}. Physica D 213 (2006) pp. 15-24.

\bibitem{P} Perdomo \emph{A small variation of the Taylor method and periodic solution of the 3-body problem}. http://arxiv.org/pdf/1410.1757.pdf

 \end{thebibliography}
\end{document}